\documentclass[12pt]{article}

\usepackage[cp1251]{inputenc}
\usepackage{amsmath,amssymb}
\usepackage[english,russian,ukrainian]{babel}

\Ukrainian

\begin{document}

\large

\begin{center}
{\LARGE Solving of partial differential equations under minimal
conditions}
\end{center}

\vskip 0.3cm

\centerline{\large Volodymyr Kyrylovych Maslyuchenko, Volodymyr
Vasylyovych Mykhaylyuk}

\vskip 0.3cm

(Chernivtsi National University, Department of Applied
Mathematics, Chair of Mathematical Analysis, Chernivtsi, vul.
Kotsyubyns'koho 2, 58012, phone (0372) 584888, e-mail:
vmykhaylyuk@ukr.net)





\vskip 0.5cm

V.K.Maslyuchenko, V.V. Mykhaylyuk {\it Solving of partial
differential equations under minimal conditions.}

It is proved that a differentiable with respect to each variable
function $f:\mathbb R^2\to\mathbb R$ is a solution of the equation
$ \frac{\partial u}{\partial x} + \frac{\partial u}{\partial y}=0$
if and only if there exists a function $\varphi:\mathbb
R\to\mathbb R$ such that $f(x,y)=\varphi(x-y)$. This gives a
positive answer to a question of R.~Baire. Besides, we use this
result to solving analogous partial differential equations in
abstract spaces and partial differential equations of
higher-order.

\vskip 0.5cm

{\bf 1. Introduction.}

\vspace{0,5cm}

Let $X, Y, Z$ be an arbitrary sets and $f:X\times Y\to Z$. For any
$x\in X$ and $y\in Y$ we define mappings $f^x:Y\to Z$ and
$f_y:X\to Z$ by the following equalities: $f^x(y)=f_y(x)=f(x,y)$.
We say that {\it a mapping $f$ separately has $P$} for some
property $P$ of mappings (continuity, differentiability, etc.) if
for any $x\in X$ and $y\in Y$ the mappings $f^x$ and $f_y$ have
$P$.

R.~Baire in fifth section of his PhD thesis [1] raised a problem
of solving of differential equations with partial derivatives
under minimal requirements, that is, a problem of solving of some
differential equation in the class of functions satisfied strongly
necessary conditions for the existence of expressions which are
contained in this equation. Besides, considering the equation
$$
\frac{\partial f}{\partial x} + \frac{\partial f}{\partial y}=0,
\eqno(1)
$$
he proved, using rather laborious arguments, that a jointly
continuous separately differentiable function $f:\mathbb
R^2\to\mathbb R$ is a solution of (1) if and only if there exists
a differentiable function $\varphi:\mathbb R\to\mathbb R$ such
that $f(x,y)=\varphi(x-y)$ for any $x,y\in \mathbb R$. Taking into
account the solution of this equation in the class of
differentiable functions $f$ (which can be obtained by introducing
of new variables $t=x-y$ and $s=x+y$), the given result means that
every jointly continuous separately differentiable solution of (1)
is differentiable. It is clear that the continuity condition on
$f$ is not necessary for the existence of partial derivatives of
$f$. Hence R.~Baire naturally raised the following question.

{\bf Question 1.1 (R.~Baire, [1, p.118]).} {\it Let $f:\mathbb
R^2\to\mathbb R$ be a separately differentiable solution of (1).
Does there exist a differentiable function $\varphi:\mathbb
R\to\mathbb R$ such that $f(x,y)=\varphi(x-y)$ for any $x,y\in
\mathbb R$?}

Note that a result analogous to Baire's result was independently
obtained in [2] where Question 1.1 was formulated too. Note that
the method used in [2] is based essentially on the joint
continuity of $f$; it is very nice and rather simpler than the
method from [1]. But, in fact, R.Baire in [1] solved (1) for
separately differentiable functions $f$ which are continuous on
every line $y=x+c$ (see Theorem 4.1).

Besides, some results concerning solutions of the following
equation
$$
\frac{\partial f}{\partial x}\cdot \frac{\partial f}{\partial y}=0
\eqno(2)
$$
appeared to the end of XX century.

So, it was proved in [3] that every continuously differentiable
solution $f:\mathbb R^2\to\mathbb R$ of equation (2) depends only
on one variable. Besides, this result was carried over the
mappings $f:X\times Y\to Z$ with locally convex range space $Z$.
Also it was shown the essentiality of the local convexity of the
space $Z$. An analogous result for separately differentiable
functions was obtained in [4]. Moreover, using rather delicate
topological arguments, it was proved that if $f:\mathbb
R^2\to\mathbb R$ is a separately continuous function and for every
point $p\in\mathbb R^2$ there exists at least one of partial
derivatives $\frac{\partial f}{\partial x}(p)$ and $\frac{\partial
f}{\partial y}(p)$, and it is equal to zero, then $f$ depends only
on one variable. This result from [4] was generalized in [5] to
the case of so-called separately $L$-differentiable mappings
$f:X\times Y\to Z$, where $X, Y,Z$ are real vector spaces and $L$
is a subspace of the space of all linear functionals on $Z$ which
separates points from $Z$.

In given paper we firstly develop a technique from [1] and study
properties of separately differentiable vector-valued functions of
two real variables (Section 2). Further, in Section 3 we establish
prove necessary and sufficient conditions under which
metric-valued functions defined on an interval are constant. Also
we obtain the following property of separately pointwise Lipschitz
(in particular, separately differentiable) functions: the
restriction of such a function on an arbitrary set has nowhere
dense discontinuity point set. This property makes possible to
give a positive answer to Question 1.1. In two last sections this
result we generalize to the case of mappings defined on the square
of a vector space and then we apply it to solve partial
differential equations of higher-orders.

\vspace{0,5cm}

{\bf 2. Auxiliary Baire function and separately differentiable
functions on $\mathbb R^2$.}

\vspace{0,5cm}

In this section we introduce an auxiliary function which connected
with the difference relation analogously as in [1] for real
functions, study its properties and use it for investigation of
separately differentiable functions.

For arbitrary $a,b \in\mathbb R$ with $a<b$ we denote by $[a;b]$,
$[a;b)$, $(a;b]$ and $(a;b)$ the corresponding intervals on
$\mathbb R$.

Let $Z$ is a vector space and $f:\mathbb R\to Z$ is a function.
For arbitrary $x,y\in \mathbb R$, $x\ne y$, and $B\subseteq Z$ put
$r_f(x,y)=\frac{f(x)-f(y)}{x-y}$ and $\Delta(B,f,x)=\{\delta\in
(0;1]:\,(\forall p',p''\in (x-\delta; x)\times
(x;x+\delta))\,\,\,(r_f(p')-r_f(p'')\in B)\}$.

Define a function $\lambda(B,f):\mathbb R\to \mathbb R$ by
following: $\lambda(B,f)(x) = {\rm sup} \Delta(B,f,x)$ if
$\Delta(B,f,x)\ne \O$ and $\lambda(B,f)(x) = 0$ if
$\Delta(B,f,x)=\O$.

Let $Z$ is a Hausdorff topological vector space. A mapping
$f:\mathbb R\to Z$ is called {\it differentiable at a point
$x_0\in\mathbb R$} if there exists $f'(x_0)=\lim\limits_{x\to
x_0}\frac{f(x)-f(x_0)}{x-x_0}$. Note that for a topological vector
space $Z$, a differentiable at $x_0$ function $f:\mathbb R\to Z$
and an arbitrary neighborhood  $W$ of zero in $Z$ we have
$\lambda(W,f)(x_0)>0$. Besides, putting $r_f(x_0,x_0)=f'(x_0)$ we
obtain that $\Delta(W,f,x_0)=\{\delta\in (0;1]:\,(\forall
p',p''\in (x_0-\delta; x_0]\times
[x_0;x_0+\delta))\,\,\,(r_f(p')-r_f(p'')\in W)\}$ for any closed
neighborhood $W$ of zero in $Z$ .

{\bf Theorem 2.1.} {\it Let $Z$ be a Hausdorff topological vector
space, $f:\mathbb R^2\to Z$ be a differentiable in the first
variable and continuous in the second variable function and $W$ be
a closed neighborhood of zero in $Z$. Then the function $g:\mathbb
R^2\to\mathbb R$, $g(x,y)=\lambda(W,f_y)(x)$, is  an jointly upper
semi-continuous function.}

P r o o f . Let $x_0,y_0\in \mathbb R$, $\gamma=g(x_0,y_0)$ and
$\varepsilon >0$. If $\gamma + \varepsilon
> 1$, then $g(x,y)\leq1<\gamma + \varepsilon$ for every $x,y\in
\mathbb R$.

Now let $\gamma + \varepsilon \leq 1$. Then $\delta_0 = \gamma +
\frac{\varepsilon}{3}\leq 1$. Since $g(x_0,y_0)<\delta_0$,
$\delta_0\not\in A(W,f_{y_0},x_0)$. Therefore there exist
$x_1,x'_1 \in (x_0-\delta_0;x_0)$ and $x_2,x'_2\in
(x_0;x_0+\delta_0)$ such that
$$
\frac{f(x_2,y_0)- f(x_1,y_0)}{x_2-x_1} - \frac{f(x'_2,y_0)-
f(x'_1,y_0)}{x'_2-x'_1}\not\in W.
$$
The continuity of $f$ in the second variable and the closedness of
$W$ imply the existence of a neighborhood $V$ of $y_0$ in $\mathbb
R$ such that
$$
\frac{f(x_2,y)- f(x_1,y)}{x_2-x_1} - \frac{f(x'_2,y)-
f(x'_1,y)}{x'_2-x'_1}\not\in W
$$
for every $y\in V$. Put
$s=\min\{x_0-x_1,x_0-x'_1,x_2-x_0,x'_2-x_0,\frac{\varepsilon}{3}\}$,
$U=(x_0-s;x_0+s)$ and $\delta_1=\gamma+\frac{2\varepsilon}{3}$.
Then $x_1,x'_1\in (x-\delta_1;x)$ and $x_2,x'_2\in (x;x+\delta_1)$
for every $x\in U$. Therefore $\delta_1\not\in A(W,f_y,x)$ and
$g(x,y)\leq \delta_1 < \gamma+\varepsilon$ for every $x\in U$ and
$y\in V$.

Thus $g$ is a jointly upper semi-continuous at $(x_0,y_0)$
function.\hfill$\diamondsuit$

Let $q,p\in \mathbb R^2$. The Euclid distance in $\mathbb R^2$
between $q$ and $p$ we denote by $d(q,p)$. If $q\ne p$ then by
$\alpha(q,p)$ we denote the angle between the vector
$\overrightarrow{pq}$ and the positive direction of abscissa.

The following theorem shows that using the function $\lambda$ one
can obtain some properties of separately differentiable functions.

{\bf Theorem 2.2.} {\it Let $Z$ be a topological vector space,
$f:\mathbb R^2\to Z$ be a separately differentiable function,
$E\subseteq \mathbb R^2$ be a nonempty set and $W$ be an arbitrary
neighborhood of zero in $Z$. Then for any open in  $E$ nonempty
set $G$ there exists point $p_0\in G$ and its neighborhood  $O$ in
$E$ such that for any distinct points $p,q\in O$ the following
inclusion holds:
$$
\frac{f(q)-f(p)}{d(q,p)} - \left(f'_x(p_0)\cos\alpha(q,p) +
f'_y(p_0)\sin\alpha(q,p)\right)\in W.
$$}

P r o o f . Note that it is sufficient to consider the case of
closed set $E$.

Let $G\subseteq E$ is an arbitrary nonempty open in $E$ set and
$W_1$ be such closed radial neighborhood of zero in $Z$ that
$W_1+W_1+W_1+W_1+W_1+W_1\subseteq W$. Consider functions
$g_1:\mathbb R^2\to \mathbb R$, $g_1(x,y)=\lambda(W_1,f_y)(x)$ and
$g_2:\mathbb R^2\to \mathbb R$, $g_2(x,y)=\lambda(W_1,f^x)(y)$.
According to Theorem 2.1, $g_1$ and $g_2$ are jointly upper
semi-continuous. For every $n\in \mathbb N$ put $E_n=\{(x,y)\in E:
g_1(x,y)\geq \frac{1}{n}, g_2(x,y)\geq \frac{1}{n}\}$. Evidently,
all the sets $E_n$ are closed in a Baire space $E$. Since
$g_1(x,y)>0$ and $g_2(x,y)>0$ for any $(x,y)\in \mathbb R^2$,
$E=\bigcup\limits_{n=1}^{\infty}E_n$. Then there exist an open in
$G$ nonempty set $H\subseteq G$ and $n_0\in \mathbb N$ such that
$H\subseteq E_{n_0}$.

Fix an arbitrary point $p_0=(x_0,y_0)\in H$. Denote
$x_1'=x_0-\frac{1}{n_0}$, $x_2'=x_0+\frac{1}{n_0}$,
$y_1'=y_0-\frac{1}{n_0}$ and $y_2'=y_0+\frac{1}{n_0}$. Separately
continuity of $f$ implies that there exists $\delta
<\frac{1}{2n_0}$ such that
$$
r_{f^{x_0}}(y_2',y_1') - r_{f^{x}}(y_2',y_1') \in
W_1\,\,\,\,\mbox{and}\,\,\,\,r_{f_{y_0}}(x_2',x_1') -
r_{f_{y}}(x_2',x_1') \in W_1
$$
for any $x\in U=(x_0-\delta;x_0+\delta)$ and $y\in
V=(y_0-\delta;y_0+\delta)$. Put $O=(U\times V)\cap H$.

Let $p=(x_1,y_1)$, $q=(x_2,y_2)$ are distinct points from the set
 $O$ and $\alpha=\alpha(q,p)$. If
$x_1\ne x_2$ and $y_1\ne y_2$ then
$$
\frac{f(q)-f(p)}{d(q,p)}=\frac{f(q)-f(x_1,y_2)}{x_2-x_1}\cdot\frac{x_2-x_1}{d(q,p)}
\, +
$$
$$
+\, \frac{f(x_1,y_2)-f(p)}{y_2-y_1}\cdot\frac{y_2-y_1}{d(q,p)} =
r_{f_{y_2}}(x_2,x_1)\cos\alpha + r_{f^{x_1}}(y_2,y_1)\sin\alpha.
$$
If  $x_1=x_2$ or $y_1=y_2$ then $\cos \alpha = 0$ or $\sin\alpha =
0$ respectively and therefore
$$
\frac{f(q)-f(p)}{d(q,p)}= r_{f_{y_2}}(x_2,x_1)\cos\alpha +
r_{f^{x_1}}(y_2,y_1)\sin\alpha.
$$

Since $p_0\in E_{n_0}$, besides, $g_1(p_0)\geq\frac{1}{n_0}$,
there exists $\delta_1>\frac{1}{2n_0}$ such that $\delta_1\in
\Delta(W_1,f_{y_0},x_0)$. Hence
$$
r_{f_{y_0}}(x_2',x_1') - r_{f_{y_0}}(x_0,x_0) \in W_1,
$$
provided $x_1'\in(x_0-\delta_1;x_0]$, $x_2'\in[x_0;x_0+\delta_1)$
and $f'_x(p_0)=r_{f_{y_0}}(x_0,x_0)$.

Note also that $q\in E_{n_0}$, besides, $g_1(q)\geq\frac{1}{n_0}$.
Since $\frac{1}{2n_0}+\delta<\frac{1}{n_0}$, there exists
$\delta_2\geq\frac{1}{2n_0}+\delta > 2\delta$ such that
$\delta_2\in \Delta(W_1,f_{y_2},x_2)$. Then
$x_1'=x_0-\frac{1}{2n_0}<x_0-\delta<x_2$, $x_2-x_1'<
x_0+\delta-x_0+\frac{1}{2n_0}\leq \delta_2$,
$x_2'=x_0+\frac{1}{2n_0}> x_0+\delta>x_2$ and
$x_2'-x_2<x_0+\frac{1}{2n_0} - x_0+\delta\leq \delta_2$. Thus
$x_1'\in (x_2-\delta_2; x_2]$ and $x_2'\in[x_2;x_2+\delta_2)$.
Inequalities $|x_1-x_2|<2\delta<\delta_2$ implies
$$
r_{f_{y_2}}(x_2,x_1) - r_{f_{y_2}}(x_2',x_1') \in W_1.
$$

Since $y_2\in V$,
$$
r_{f_{y_2}}(x_2',x_1') - r_{f_{y_0}}(x_2',x_1') \in W_1.
$$

Now we have
$$
r_{f_{y_2}}(x_2,x_1) - f'_x(p_0) = \left(r_{f_{y_2}}(x_2,x_1) -
r_{f_{y_2}}(x_2',x_1')\right) +
$$
$$
\left(r_{f_{y_2}}(x_2',x_1') -
 r_{f_{y_0}}(x_2',x_1')\right) + \left(r_{f_{y_0}}(x_2',x_1') -
f'_x(p_0)\right)\in W_1+W_1+W_1.
$$
Analogously
$$
r_{f^{x_1}}(y_2,y_1) - f'_y(p_0) \in W_1+W_1+W_1.
$$
Then
$$
\frac{f(q)-f(p)}{d(q,p)} - (f'_x(p_0)\cos\alpha +
f'_y(p_0)\sin\alpha) =
$$
$$
= \cos\alpha\left( r_{f_{y_2}}(x_2,x_1) - f'_x(p_0)\right) +
\sin\alpha \left(r_{f^{x_1}}(y_2,y_1) - f'_y(p_0)\right) \in
$$
$$
\in \cos\alpha (W_1+W_1+W_1) + \sin\alpha (W_1+W_1+W_1)\subseteq
W.
$$

This complete the proof. \hfill$\diamondsuit$

\vspace{0,5cm}

{\bf 3. Separately pointwise Lipschitz functions and pointwise
changeable functions.}

\vspace{0,5cm}

Firstly recall some definitions.

Let $(X, |\cdot - \cdot|_X)$ and $(Y, |\cdot - \cdot|_Y)$ be
metric spaces. A mapping $f:X\to Y$ {\it satisfies Lipschitz
condition with a constant $C>0$} if $|f(x) - f(y)|_Y\leq C
|x-y|_X$ for any $x,y\in X$. A mapping $f:X\to Y$ is called {\it
pointwise Lipschitz} if for any point $x_0\in X$ there exist a
neighborhood $U$ of point $x_0$ in $X$ and $C>0$ such that
$|f(x_0) - f(x)|_Y\leq C |x_0-x|_X$ for any $x\in U$. A mapping
$f:X\to Y$ is called {\it pointwise changeable}, if for every
$\varepsilon>0$ the union $G_{\varepsilon}$ of the system ${\cal
G}_{\varepsilon}$ of all open nonempty sets $G\subseteq X$ such
that $f|_G$ satisfies the Lipschitz property with the constant
$\varepsilon$, is an everywhere dense set.

The following property of separately pointwise Lipschitz mappings
plays an important role in obtaining of a positive answer to
Question 1.1.

{\bf Theorem 3.1.} {\it Let $(X, |\cdot - \cdot|_X)$ and $(Y,
|\cdot - \cdot|_Y)$ be metric spaces such that the space $X\times
Y$ is a hereditarily Baire space, $(Z, |\cdot - \cdot|_Z)$ be a
metric space and $f:X\times Y\to Z$ be a separately pointwise
Lipschitz mapping. Then for any nonempty set $E\subseteq X\times
Y$ the discontinuity point set $D(f|_E)$ of mapping $f|_E$ is
nowhere dense in $E$.}

P r o o f . Note that it is sufficient to prove the theorem for
closed set $E$.

Let $E\subseteq X\times Y$ be a closed nonempty set and
$G\subseteq X\times Y$ be an open set such that $W_0=G\cap
E\ne\O$. For any $n,m\in\mathbb N$ denote by $E_{nm}$ the set of
all points $(x,y)\in W_0$ such that
$$
|f(x',y)-f(x,y)|_Z\leq n|x'-x|_X {\,\,\,and\,\,\,}
|f(x,y')-f(x,y)|_Z\leq n|y'-y|_Y
$$
for any $x'\in X$ with $|x'-x|_X<\frac{1}{m}$ and $y'\in Y$ with
$|y'-y|_Y<\frac{1}{m}$. Since $f$ is separately pointwise
Lipschitz function,  $W_0=\bigcup\limits_{n,m=1}^{\infty} E_{nm}$.
We obtain that there exist $n_0,m_0\in\mathbb N$ and an open in
$E$ nonempty set $W\subseteq W_0$ such that $E_{n_0m_0}$ is dense
in $W$, provided $W_0$ is open set in a Baire space $E$.

Choose open balls $U_1$ and $V_1$ with radius $\frac{1}{2m_0}$ in
spaces  $X$ and $Y$ respectively, such that  $W_1=(U_1\times
V_1)\cap W\ne\O$. Let's show that function $f$ satisfies Lipschitz
condition on the set  $W_1$ with the constant $2n_0$ with respect
to the maximum-metric $|\cdot-\cdot|_{X\times Y}$ on $X\times Y$.

Let $p_1=(x_1,y_1), p_2=(x_2,y_2)\in W_1$. Fix arbitrary
$\varepsilon, \delta >0$. Since $f$ is continuous in the first
variable at points $p_1$ and $p_2$ and the set $E_{n_0m_0}$ is
dense in $W_1$, there exist
$(\tilde{x}_1,\tilde{y}_1),(\tilde{x}_2,\tilde{y}_2)\in W_1\cap
E_{n_0 m_0}$ such that
$$
|x_1-\tilde{x}_1|_X<\delta,\,\,|y_1-\tilde{y}_1|_Y<\delta,\,\,|x_2-\tilde{x}_2|_X<\delta,
\,\,|y_2-\tilde{y}_2|_Y<\delta,
$$
$$
|f(x_1,y_1)-f(\tilde{x}_1,y_1)|_Z<\varepsilon \,\,\,{\mbox
and}\,\,\,|f(x_2,y_2)-f(\tilde{x}_2,y_2)|_Z<\varepsilon.
$$
Then
$$
|f(p_1)-f(p_2)|_Z \leq |f(x_1,y_1)-f(\tilde{x}_1,y_1)|_Z +
|f(\tilde{x}_1,y_1)-f(\tilde{x}_1,\tilde{y_1})|_Z +
$$
$$
|f(\tilde{x}_1,\tilde{y}_1)-f(\tilde{x}_1,\tilde{y}_2)|_Z +
|f(\tilde{x}_1,\tilde{y}_2)-f(\tilde{x}_2,\tilde{y}_2)|_Z +
|f(\tilde{x}_2,\tilde{y}_2)-f(\tilde{x}_2,y_2)|_Z +
$$
$$
|f(\tilde{x}_2,y_2)-f(x_2,y_2)|_Z \leq \varepsilon + n_0
|y_1-\tilde{y}_1|_Y + n_0|\tilde{y}_1-\tilde{y}_2|_Y +
n_0|\tilde{x}_1-\tilde{x}_2|_X +
$$
$$
n_0 |\tilde{y}_2 - y_2|_Y + \varepsilon \leq 2\varepsilon +
2\delta n_0 + n_0(|y_1-y_2|_Y + 2\delta) + n_0(|x_1-x_2|_X +
2\delta) =
$$
$$
2\varepsilon + 6\delta n_0 + n_0(|x_1-x_2|_X + |y_1-y_2|_Y)\leq
2\varepsilon + 6\delta n_0 + 2n_0|p_1-p_2|_{X\times Y}.
$$
Tending $\varepsilon$ and $\delta$ to zero, we obtain
$$
|f(p_1) - f(p_2)|_Z \leq 2n_0 |p_1 - p_2|_{X\times Y}.
$$

Hence, $f|_E$ is continuous on the set $W_1$.\hfill$\diamondsuit$

Note that obtained property of separately pointwise Lipschitz
mappings is new, but for real-valued separately differentiable
functions of two variables this property can be obtained from the
analog of Theorem 2.2, which was presented in [1]. Besides, in [6]
it was proved that the discontinuity point set of function of two
real variables, which is differentiable in the first variable and
continuous in the second one, is nowhere dense. This result was
generalized in [7].

For a topological space $X$ and a set $A\subseteq X$ by
$\overline{A}$ we denote closure of $A$ in $X$.

The following characterization was obtained in [1] for real-valued
functions of one real variable.

{\bf Theorem 3.2.} {\it Let $X\subseteq\mathbb R$ be a nonempty
interval, \mbox{$(Y,|\cdot - \cdot|_Y)$} be a metric space,
$f:X\to Y$ be a continuous pointwise changeable on every closed
set mapping. Then $f$ is constant.}

P r o o f . For any $x_1,x_2\in X$, $x_1\ne x_2$ put $r(x_1,x_2)=
\frac{|f(x_2) - f(x_1)|_Y}{|x_2 - x_1|}$ and for every $x\in X$
put $\displaystyle g(x) =\inf\limits_{\delta>0}\sup
\{r(x_1,x_2):x-\delta<x_1<x_2<x+\delta\}$.

Let's show that for any $a,b\in X$, $a<b$, there exists point
$c\in [a;b]$ such that $g(c)\geq r(a,b)$.

Let $a\leq x < y < z \leq b$. Then $r(x,z)\leq
\frac{y-x}{z-x}r(x,y) + \frac{z-y}{z-x}r(y,z)$. Hence, $r(x,z)\leq
r(x,y)$ or $r(x,z)\leq r(y,z)$. Now it is easy to construct such
sequence $(I_n)_{n=1}^{\infty}$ of segments
$I_n=[a_n;b_n]\subseteq [a;b]$ such that
$\lim\limits_{n\to\infty}(b_n-a_n)=0$, $I_{n+1}\subseteq I_n$ and
$r(a_n,b_n)\geq r(a,b)$ for every $n\in\mathbb N$. Then for point
$c\in \bigcap\limits_{n=1}^{\infty} I_n$ we have $g(c)\geq
r(a,b)$.

Therefore, if $a,b\in X$, $\varepsilon>0$ and $g(x)\leq
\varepsilon$ for any  $x\in (a,b)\subseteq X$, then
$r(x,y)\leq\varepsilon$ for any $x,y\in (a;b)$, what implies
$r(a,b)\leq \varepsilon$, provided  $f$ is continuous.

Assume that $f$ is not constant. Then there exists $\varepsilon
>0$ such that $\displaystyle E=\{x\in X: g(x)>\varepsilon\}\ne\O$.
Since $f$ is pointwise changeable on the set $F=\overline E$,
there exist $x_0\in E$ and $\delta>0$ such that
$r(x,y)<\varepsilon$ for any distinct $x,y\in F\cap U$, where
$U=(x_0-\delta; x_0+\delta)$.

Let $x,y\in U$ be arbitrary distinct points. Let's show that
 $r(x,y)\leq\varepsilon$. Let us assume that
$x<y$. Firstly consider the case of $x,y\not\in F$. If $(x;y)\cap
F=\O$, then $(x;y)\cap E=\O$ and $r(x,y)\leq\varepsilon$. Let
$(x;y)\cap F\ne\O$. Choose points $u,v\in F$ such that $x< u \leq
v < y$, $(x;u)\cap F=\O$ and $(v;y)\cap F=\O$. Then, as in above,
$r(x,u)\leq\varepsilon$ and $r(v,y)\leq\varepsilon$. If $u<v$,
then
$$
r(x,y)=\frac{u-x}{y-x}r(x,u)+\frac{v-u}{y-x}r(u,v)+\frac{y-v}{y-x}r(v,y),
$$
therefore $r(x,y)\le\varepsilon$. When $u=v$ we use the equality

$$
r(x,y)=\frac{u-x}{y-x}r(x,u)+\frac{y-u}{y-x}r(u,y).
$$

In the case of $x\in F$ or $y\in F$ we use analogous reasons.

Thus, $ \displaystyle {\rm
sup}\{r(x,y):x_0-\delta<x<y<x_0+\delta\}\le\varepsilon$. Then
$g(x_0)\le\varepsilon$, what contradicts to $x_0\in
E$.\hfill$\diamondsuit$

{\bf Corollary 3.3.} {\it Let $X$ be an arbitrary normed space,
$(Y,|\cdot - \cdot|_Y)$ be a metric space, $f:X\to Y$ be a
continuous pointwise changeable on every closed set mapping. Then
$f$ is constant.}

P r o o f . It is enough to prove that $f(x)=f(0)$ for any $x\in
X$.

Let $x_0\in X$, $x_0\ne 0$, is an arbitrary point. Consider the
function $g:\mathbb R\to Y$, $g(\alpha)= f(\alpha x_0)$. Since
$|\alpha - \beta|= \frac{1}{\|x_0\|}\|\alpha x_0 - \beta x_0\|$
and $f$ is continuous pointwise changeable on every closed set
mapping,  $g$ satisfies conditions of Theorem 3.2. Therefore, $g$
is constant and, besides, $f(x_0)= g(1) = g(0) =
f(0)$.\hfill$\diamondsuit$

The following two examples demonstrate that there is no analogous
property for mappings defined on an arbitrary metric space, and
from the other side, this property has no any equivalent
formulation in topological terms.

E x a m p l e  3.4 .  Let $(X,|\cdot - \cdot|_X)$ be a metric
space with the discrete metric, i.e. $|x_1-x_2|_X=1$ when $x_1\ne
x_2$, and $(Y, |\cdot - \cdot|_Y)$ be an arbitrary metric space.
Then every mapping $f:X\to Y$ is continuous and pointwise
changeable on every closed set.

E x a m p l e   3.5 . Let $0<p<1$ and $\mathbb R_p$ be the real
line with the metric $|x-y|_p=|x-y|^p$. Then the identical map
$f:\mathbb R_p\to\mathbb R$, $f(x)=x$, is an pointwise changeable
on every closed set homeomorphism.

\vspace{0,5cm}

{\bf 4. The equation $f'_x + f'_y =0$.}

\vspace{0,5cm}

In this section we give a positive answer to Question 1.1.

Actually, the following theorem was proved in [1], but R.~Baire
instead of continuity of function $f$ on respective lines put on
$f$ stronger condition of joint continuity.

{\bf Theorem 4.1.} {\it Let $f:\mathbb R^2\to\mathbb R$ be a
separately differentiable function and $c\in \mathbb R$ such that
the restriction of function $f$ to the set $A=\{(x,y)\in\mathbb
R^2: y-x=c\}$ is continuous and $f'_x(p)+f'_y(p)=0$ for every
$p\in A$. Then the function $g:\mathbb R\to\mathbb R$,
$g(x)=f(x,c+x)$, is constant.}

P r o o f . Since $\cos\alpha(q,p)=\sin\alpha(q,p)$ for any
distinct points $p,q\in A$, Theorem 2.2 implies that the
continuous function $g$ is pointwise changeable on every closed
set. It remains to apply Theorem 3.2. \hfill$\diamondsuit$

In the proof of the main result we will use the following
auxiliary fact.

{\bf Lemma 4.2.} {\it Let $I=(a;b)\subseteq \mathbb R$ be an
arbitrary nonempty interval, $c\in \mathbb R$, $\delta
>0$, $W=\{(x,y)\in \mathbb R^2: x\in I, |y-x-c|\leq\delta\}$,
$f:\mathbb R^2\to\mathbb R$ and $g:\mathbb R^2\to \mathbb R$ be
such separately continuous functions that $f(x,y)=g(x,y)$ for any
$(x,y)\in W$. Then $f(x,y)=g(x,y)$ for any $(x,y)\in\overline W$.}

P r o o f . Let $x_0=a$ and $|y_0-x_0-c|<\delta$. Then $f(x_0,y_0)
= \lim\limits_{x\to a+0} f(x,y_0) = \lim\limits_{x\to a+0}
g(x,y_0) = g(x_0,y_0)$. Analogously, if $x_0=b$ and
$|y_0-x_0-c|<\delta$, then $f(x_0,y_0)=g(x_0,y_0)$.

Now let $x_0=a$ and $y_0-x_0-c=\delta$. Then
$f(x_0,y_0)=\lim\limits_{y\to y_0-\,0}f(x_0,y)=\lim\limits_{y\to
y_0-\,0}g(x_0,y)=g(x_0,y_0)$. We use analogous reasons in the case
of $x_0=a$ and $y_0-x_0-c=-\delta$, or $x_0=b$ and
$y_0-x_0-c=\pm\delta$. \hfill$\diamondsuit$

Let  $X$ be a topological space, $x_0\in X$, ${\cal U}$ be a
system of all neighborhoods of point $x_0$ in $X$, $(Y,|\cdot -
\cdot|_Y)$ be a metric space and $f:X\to Y$. Recall that a real
$\omega_f(x_0) = \inf\limits_{U\in {\cal U}}\sup\limits_{x',x''\in
U}|f(x')-f(x'')|_Y$ is called  {\it the oscillation of mapping $f$
at $x_0$}.

Now let us proof our main result.

{\bf Theorem 4.3.} {\it Let $f:\mathbb R^2\to\mathbb R$ be a
separately differentiable function such that  $f'_x(p)+f'_y(p)=0$
for every $p\in \mathbb R^2$. Then for any $c\in \mathbb R$ the
function $f$ is constant on the set $A=\{(x,y)\in\mathbb R^2:
y-x=c\}$.}

P r o o f .  According to Theorem  4.1 it is enough to prove that
 $f$ is continuous.

Assume that the discontinuity point set $E$ of function $f$ is
nonempty. Theorem 3.1 implies that there exists point
$p_0=(x_0,y_0)\in E$, in which the function  $f|_E$ is continuous.
Denote $\varepsilon=\omega_f(p_0)$, $c_0=y_0-x_0$ and choose
$\delta_1, \delta_2>0$ such that for any point  $p\in E\bigcap W$,
where $W=\{(x,y)\in
\mathbb{R}^2:\,\,|x-x_0|<\delta_1,\,\,|y-x-c_0|<\delta_2\}$, the
inequality $|f(p)-f(p_0)|\leq \frac{\varepsilon}{3}$ holds. Note
that for any point $q\in W$ with $|f(q)-f(p_0)|>
\frac{\varepsilon}{3}$ the function $f$ is continuous at $q$.

Consider the continuous function $g:\mathbb{R}^2 \to\mathbb{R}$,
$g(x,y)=f(x_0,x_0+y-x)$ and show that $f(p)=g(p)$ for any point
 $p\in W$ with
$|f(p)-f(p_0)|>\frac{\varepsilon}{3}$.

Let $p_1=(x_1,y_1)\in W$, besides,
$|f(p_1)-f(p_0)|>\frac{\varepsilon}{3}$. Choose $\delta>0$ such
that $|f(x_1,y)-f(x_0,y_0)|>\frac{\varepsilon}{3}$ and $(x_1,y)\in
W$ for any $y\in [y_1-\delta;y_1+\delta]$.

Denote by $\cal I$ the system of all nonempty open intervals
$I\subseteq (x_0-\delta_1;x_0+\delta_1)$ such that $x_1\in I$ and
$|f(x,y)-f(x_0,y_0)|>\frac{\varepsilon}{3}$ for any $x\in I$ and
$y\in \mathbb R$ with $|y-x-c_1|\leq\delta$, where $c_1=y_1-x_1$.
Note that $f$ is continuous at every point of compact set
$K=\{(x_1,y):y\in[y_1-\delta;y_1+\delta]\}$. Therefore the system
$\cal I$ is nonempty.

Put $I_0=(a;b)=\bigcup\limits_{I\in{\cal I}}I$ and $W_1=\{(x,y)\in
\mathbb R^2: x\in I_0, |y-x-c_1|\leq\delta\}$. Since
$|f(p)-f(p_0)|>\frac{\varepsilon}{3}$ for every $p\in W_1\subseteq
W$, the function $f$ is continuous at every point from $W_1$.
According to Theorem 2.2 the function $\varphi(x)=f(x,c+x)$ is
pointwise changeable on  $I$, and therefore, accordingly to
Theorem 3.2, $\varphi$ is constant on $I$ for every $c\in
[c_1-\delta; c_1+\delta]$, i.e. $f(x,y)= f(x,x+y-x) =
f(x_1,x_1+y-x)$ for any $(x,y)\in W_1$.

Let us show that $I_0=(x_0-\delta_1;x_0+\delta_1)$. Assume that
$a>x_0-\delta_1$. Then Lemma 4.2 implies $f(x,y)=f(x_1,x_1+y-x)$
for any $(x,y)\in \overline{W_1}$, besides,
$f(a,y)=f(x_1,x_1+y-a)$ for any $y\in
[a+c_1-\delta;a+c_1+\delta]$. Note that $(x_1;x_1+y-a)\in K$ if
$y\in [a+c_1-\delta;a+c_1+\delta]$, then $|f(p)-f(p_0)|>
\frac{\varepsilon}{3}$ for every $p\in K_1$, where $K_1=\{(a,y):
|y-a-c_1|\leq\delta\}$. Since $K_1\subseteq W$, the function $f$
is continuous at every point from the set $K_1$. Hence, there
exists a nonempty interval $I_1\subseteq
(x_0-\delta_1;x_0+\delta_1)$ such that  $a\in I_1$ and
$|f(x,y)-f(x_0,y_0)|>\frac{\varepsilon}{3}$ for any $x\in I_1$ and
$y\in \mathbb R$ with $|y-x-c_1|\leq\delta$. Then  $I_0\bigcup
I_1\in{\cal I}$, what contradicts to the definition of the set
$I_0$. We use analogous reasons if $b<x_0+\delta$.

Thus, $I_0=(x_0-\delta_1;x_0+\delta_1)$. Then $(x_0,x_0+c_1)\in
W_1$ and $f(x_0,x_0+c_1)=f(x_1,x_1+c_1) = f(x_1,y_1)$, then
$g(x_1,y_1)= f(x_0,x_0+c_1)=f(x_1,y_1)$.

Since $\omega_f(p_0)=\varepsilon$, then there exists a sequence
 $(q_n)_{n=1}^\infty$ of points
$q_n=(u_n,v_n)\in W$ such that $|f(q_n)-f(p_0)|>
\frac{\varepsilon}{3}$ and $\lim\limits_{n\to\infty}q_n=p_0$.
Then, using continuity of $g$, we obtain
$\lim\limits_{n\to\infty}f(q_n) = \lim\limits_{n\to\infty}g(q_n) =
g(p_0) = f(p_0)$. But the last equalities contradict to the choice
of $(q_n)_{n=1}^\infty$. \hfill$\diamondsuit$

{\bf Corollary 4.4.} {\it Let $k\in \mathbb R$, $k\ne 0$,
$f:\mathbb R^2\to\mathbb R$ be such a separately differentiable
function that $f'_x(p)+ kf'_y(p)=0$ for every  $p\in \mathbb R^2$.
Then there exists a differentiable function $\varphi:\mathbb
R\to\mathbb R$ such that $f(x,y)=\varphi(kx-y)$ for any  $x,y\in
\mathbb R$.}

\vspace{1,0cm}

{\bf 5. Equations for separately $L$-differentiable functions.}

\vspace{0,5cm}

In this section we apply Theorem 4.3 to solving of differential
equations in abstract spaces.

Let $X$ be a vector space, $Z$ be a set and $L$ be a system of
functions $l:Z\to \mathbb R$. We say that a mapping $f:X\to Z$ is
{\it $L$-differentiable at $x_0\in X$} if for arbitrary $h\in X$
and $l\in L$ the function $g:\mathbb R\to \mathbb R$,
$g(t)=l(f(x_0+th))$, is differentiable at $t_0=0$, i.e. there
exists $A(h,l)= \lim\limits_{t\to 0}\frac{l(f(x_0+th))-
l(f(x_0))}{t}$. The mapping $A:X\times L\to\mathbb R$ is called
{\it $L$-derivative of $f$ at $x_0$} and is denoted by $Df(x_0)$.
Besibes, $Df(x_0)(h,l)$ we denote by $Df(x_0,h,l)$.

A mapping $f:X\to Z$ is called {\it $L$-differentiable} if $f$ is
$L$-differentiable at every point $x\in X$.

Recall that a system $L$ of functions defined on a set $Z$, {\it
separates points from $Z$} if for arbitrary distinct points
$z_1,z_2\in Z$ there exists $l\in L$ such that $l(z_1)\ne l(z_2)$.

{\bf Theorem 5.1.} {\it Let $X$ be a vector space, $Z$ be a set,
$L$ be a system of functions defined on $Z$ which separates points
from $Z$ and $f:X^2\to Z$ be a separately $L$-differentiable
mappings such that
$$
Df^x(y) + Df_y(x) =0
$$
for every $x,y\in X$. Then there exists an $L$-differentiable
mapping $\varphi:X\to Z$ such that $f(x,y) = \varphi(x-y)$ for
every $x,y\in X$.}

P r o o f . Firstly show that $f(x,y) = \varphi(x-y)$ for some
mapping $\varphi:X\to Z$. It is enough to prove that $f(x_1,y_1) =
f(x_2,y_2)$ if $x_1-y_1=x_2-y_2$.

Suppose that there exist $x_1,y_1,x_2,y_2\in X$ such that
$x_1-y_1=x_2-y_2$ and $f(x_1,y_1) = z_1\ne f(x_2,y_2)=z_2$. Since
the system $L$ separates points from $Z$, there exists $l\in L$
such that $l(z_1)\ne l(z_2)$. Put $h=x_2-x_1 = y_2-y_1$ and
consider the function $u: \mathbb R^2\to\mathbb R$, $u(s,t)=
l(f(x_1+sh,y_1+th))$.

Show that $u$ is a separately differentiable function with $u'_s +
u'_t =0$.

Let $s_0,t_0\in \mathbb R$, $x_0=x_1+s_0h$ and $y_0=y_1+t_0h$.
Then
$$
u'_s(s_0,t_0) = \lim\limits_{s\to s_0}\frac{l(f(x_1+sh,y_1+t_0h))-
l(f(x_1+s_0h,y_1+t_0h))}{s-s_0} =
$$
$$
= \lim\limits_{s\to s_0}\frac{l(f(x_0+(s-s_0)h,y_0))-
l(f(x_0,y_0))}{s-s_0} =
$$
$$
= \lim\limits_{s\to s_0}\frac{l(f_{y_0}(x_0+(s-s_0)h))-
l(f_{y_0}(x_0))}{s-s_0} = Df_{y_0}(x_0,h,l).
$$

Analogously $u'_t(s_0,t_0) = Df^{x_0}(y_0,h,l)$. Since
$Df^{x_0}(y_0) + Df_{y_0}(x_0)=0$, $u'_s(s_0,t_0) + u'_t(s_0,t_0)
=0$.

Thus $u$ satisfies the conditions of Theorem 4.3 therefore
$l(z_1)=u(0,0)=u(1,1)=l(z_2)$ what contradicts to our assumption.

$L$-differentiability of $\varphi$ follows from
$\varphi(x)=f(x,0)$ and $L$-differentiability of $f_y$ if
$y=0$.\hfill$\diamondsuit$

Let $X$, $Z$ be topological vector spaces. A mapping $f:X\to Z$ is
called {\it Gateaux differentiable at point $x_0\in X$} if there
exists an linear continuous operator  $A:X\to Z$ such that
$$
\lim\limits_{t\to 0}\frac{f(x_0+th) - f(x_0)}{t} = (Ax_0)(h)
$$
for every $h\in X$. Such operator $A$ is called {\it Gateaux
derivative of mapping $f$ at point $x_0$}.

Note that for a Hausdorff topological vector space  $Z$ Gateaux
derivative is unique. A mapping $f:X\to Z$, which is Gateaux
differentiable at every point $x\in X$, is called  {\it Gateaux
differentiable}. A mapping $D$, which every Gateaux differentiable
mapping $f:X\to Z$ assigns the Gateaux derivative mapping, i.e.
$Df(x)$ is the Gateaux derivative of $f$ at point $x\in X$, we
call {\it Gateaux differentiation operator}.

{\bf Corollary 5.2.} {\it Let $X$ be a topological vector space,
$Z$ be a topological vector space such that the conjugate space
$Z^*$ separates points from $Z$ and $f:X^2\to Z$ be a mapping such
that
$$
Df^x(y) + Df_y(x) = 0
$$
for every $x,y\in X$, where $D$ is the Gateaux differentiation
operator for mapping which acting from $X$ to $Z$. Then there
exists a Gateaux differentiable mapping $\varphi:X\to Z$ such that
$f(x,y)= \varphi(x-y)$ for every $x,y\in X$.}

P r o o f . Since $Z^*$ separates points from $Z$, $Z$ is a
Hausdorff space and the definition of $D$ is correct. Besides, for
arbitrary $x,y,h\in X$ and $z^*\in Z^*$ we have
$$
\lim\limits_{t\to 0}\frac{z^*(f_y(x+th))- z^*(f_y(x))}{t} =
\lim\limits_{t\to 0}z^*(\frac{f_y(x+th)- f_y(x)}{t}) =
z^*(Df_y(x)(h))
$$
and
$$ \lim\limits_{t\to 0}\frac{z^*(f^x(y+th))- z^*(f^x(y))}{t} =
z^*(Df^x(y)(h)) = - z^*(Df_y(x)(h)).
$$
Therefore $f$ is a separately $Z^*$-differentiable mapping and
$\tilde{D}f^x(y) + \tilde{D}f_y(x) = 0$, where $\tilde{D}$ is the
$Z^*$-differentiation operator. Theorem 5.1 implies that there
exists a mapping $\varphi:X\to Z$ such that $f(x,y)=\varphi(x-y)$
for every $x,y\in X$. Since $\varphi(x)=f(x,0)$,  $\varphi$ is a
Gateaux differentiable mapping.\hfill$\diamondsuit$

\vspace{0,5cm}

{\bf 6. Higher-order equations.}

\vspace{0,5cm}

Finally we give applications of Theorem 4.3 to solving of
higher-order partial differential equations.

Let $n\in \mathbb N$ and a function $f:\mathbb R^2\to\mathbb R$
has all $n$-order partial derivatives. The sum of all $n$-order
partial derivatives of $f$ we denote by $D_nf$. Clearly that
$D_{n+k}f=D_n(D_kf)$ for every $n,k\in\mathbb N$ and any function
$f:\mathbb R^2\to\mathbb R$ which has all $(n+k)$-order partial
derivatives.

{\bf Theorem 6.1.} {\it Let $n\in \mathbb N$, a function
$f:\mathbb R^2\to\mathbb R$ has all $n$-order partial derivatives
and $D_nf(p)=0$ for every $p\in \mathbb R^2$. Then there exist
differentiable functions $\varphi_1, \dots ,\varphi_n:\mathbb
R\to\mathbb R$ such that
$$
f(x,y)=\varphi_1(x-y) + (x+y)\varphi_2(x-y) + \dots +
(x+y)^{n-1}\varphi_n(x-y)
$$
for every $x,y\in\mathbb R$.}

P r o o f . The proof is by induction of $n$.

By $n=1$ it follows from Theorem 4.3.

Assume that our assertion is true for some $n=k$ and prove it for
$n=k+1$.

Let $f:\mathbb R^2\to\mathbb R$ be a function which has all
$(k+1)$-order partial derivatives and $D_{k+1}f(p)=0$ for every
$p\in \mathbb R^2$.

Put $g=D_1f$. Since $D_kg=D_{k+1}f=0$, the assumption implies that
there exist differentiable functions $\psi_1, \dots
,\psi_k:\mathbb R\to\mathbb R$ such that
$$
g(x,y)=\psi_1(x-y) + (x+y)\psi_2(x-y) + \dots +
(x+y)^{k-1}\psi_k(x-y) .
$$

Denote $\varphi_{i+1}=\frac{1}{2i}\psi_i$ if $1\leq i \leq k$, and
put
$$
u(x,y)=(x+y)\varphi_2(x-y) + (x+y)^2\varphi_3(x-y) + \dots +
(x+y)^{k}\varphi_{k+1}(x-y).
$$
Then
$$
D_1u(x,y)= 2\varphi_2(x-y) + 4(x+y)\varphi_3(x-y) + \dots +
2k(x+y)^{k-1}\varphi_{k+1}(x-y) =
$$
$$
= \psi_1(x-y) + (x+y)\psi_2(x-y) + \dots + (x+y)^{k-1}\psi_k(x-y)
= g(x,y) = D_1f(x,y).
$$

Thus $D_1(f-u)=0$ and Theorem 4.3 implies that there exists a
differentiable function $\varphi_1:\mathbb R\to\mathbb R$ such
that $f(x,y)=\varphi_1(x,y) + u(x,y)$.\hfill$\diamondsuit$

{\bf Theorem 6.2.} {\it Let a function $f:\mathbb R^2\to\mathbb R$
has all second-order partial derivatives and
$$
f''_{xx}(p)=f''_{yy}(p)\,\,\,\,\,\mbox{and}\,\,\,\,\,f''_{xy}(p)=f''_{yx}(p)
$$
for every $p\in \mathbb R^2$. Then there exist twice
differentiable functions $\varphi, \psi:\mathbb R\to\mathbb R$
such that
$$
f(x,y)=\varphi(x+y) + \psi(x-y)
$$
for every $x,y\in\mathbb R$.}

P r o o f . Consider the function $g:\mathbb R^2\to\mathbb R$,
$g(p)=f'_x(p) - f'_y(p)$. Then $g'_x(p) + g'_y(p) = f''_{xx}(p)-
f''_{yx}(p) + f''_{xy}(p) - f''_{yy}(p) = 0$. Theorem 4.3 implies
that there exists a differentiable function $\tilde{\psi}:\mathbb
R\to \mathbb R$ such that $g(x,y)=\tilde{\psi}(x-y)$.

Pick some twice differentiable function $\psi:\mathbb R\to\mathbb
R$ such that $2\psi'=\tilde{\psi}$ and consider the function
$\tilde{f}(x,y)=f(x,y)-\psi(x-y)$. Then
$$
\tilde{f}'_x(x,y) - \tilde{f}'_y(x,y) = f'_x(x,y) - f'_y(x,y) -
2\psi'(x-y) = g(x,y) - g(x,y) = 0.
$$
Therefore Corollary 4.4 implies that there exists a differentiable
function $\varphi:\mathbb R\to \mathbb R$ such that
$\tilde{f}(x,y)=\varphi(x+y)$, i.e. $f(x,y)=\varphi(x+y) +
\psi(x-y)$ for every $x,y\in\mathbb R$. Since
$\varphi(x)=f(x,0)-\psi(x)$, $\varphi$ is a twice differentiable
function.\hfill$\diamondsuit$

R e m a r k s . The existence of $f''_{xx}$ and $f''_{yy}$ for a
function $f:\mathbb R^2\to \mathbb R$ does not imply the existence
of $f''_{xy}$ and $f''_{yx}$. For example, the Schwartz function
$$
f(x,y) = \left \{\begin{array}{rr}
 \frac{2xy}{x^2+y^2},
&
 {\rm if}\quad x^2+y^2\ne 0;
\\
  0,
&
 {\rm if}\quad x=y=0,
  \end{array} \right .
$$
is a separately infinite differentiable function, but
$f''_{xy}(0,0)$ and $f''_{yx}(0,0)$  do not exist.

In other hand, the existence of all second-order partial
derivatives of a function $f:\mathbb R^2\to \mathbb R$ and the
equality $f''_{xy}=f''_{yx}$ do not imply the jointly continuity
of $f$. Really, the function
$$
f(x,y) = \left \{\begin{array}{rr}
 \frac{2x^3y^3}{x^6+y^6},
&
 {\rm if}\quad x^6+y^6\ne 0;
\\
  0,
&
 {\rm if}\quad x=y=0,
  \end{array} \right .
$$
has all second-order partial derivatives, besides,
$f''_{xy}=f''_{yx}$ on $\mathbb R^2$ and $f$ is jointly
discontinuous at $(0,0)$.

In this connection the following question arises naturally.

{\bf Question 6.3.} {\it Let a function $f:\mathbb R^2\to\mathbb
R$ has partial derivatives $f''_{xx}$ and $f''_{yy}$ and
$$
f''_{xx}(p)=f''_{yy}(p)
$$
for every $p\in \mathbb R^2$. Do there exist twice differentiable
functions $\varphi, \psi:\mathbb R\to\mathbb R$ such that
$$
f(x,y)=\varphi(x+y) + \psi(x-y)
$$
for every $x,y\in\mathbb R$?}

\vspace{0,5cm} \normalsize

\begin{center}
{\LARGE References}
\end{center}

\begin{enumerate}

\item {\it R.Baire,} Sur les fonctions de variables reelles. -
Annali di mat. pura ed appl., ser.3., {\bf 3} (1899), 1-123.

\item {\it P.R.Chernoff, H.F.Royden} The Equation $\frac{\partial
f}{\partial x} = \frac{\partial f}{\partial y}$. - The American
Mathematical Monthly, {\bf V.82, №5} (1975), 530-531.

\item {\it V.K.Maslyuchenko,} One property of partial derivatives.
- Ukr. mat. zhurn. {\bf 39,№4} (1987), 529-531(in Russian).

\item {\it A.M.Bruckner, G.Petruska, O.Preiss, B.S.Thomson,} The
equation $u_xu_y=0$ factors. - Acta Math. Hung. {\bf 57,№3-4}
(1991), 275-278.

\item {\it A.K.Kalancha, V.K.Maslyuchenko,} A generalization of
Bruckner-Petruska-Preiss-Thomson theorem. - Mat. Studiji. {\bf
1,№1} (1999), 48-52(in Ukrainian).

\item {\it K.B$\ddot{o}$gel,} ${\rm \ddot{U}}$ber partiell
differenzierbare Funktionen. - Math. Z. {\bf 25} (1926), 490-498.

\item {\it V.H.Herasymchuk, V.K.Maslyuchenko, V.V.Mykhaylyuk,}
Varieties of Lipschitz condition and discontinuity points sets of
separately differentiable functions. - Nauk. Visn. Cherniv. Univ.
Vyp. {\bf 134}. Matematyka. Chernivtsi: Ruta, (2002), 22-29(in
Ukrainian).

\end{enumerate}

\end{document}